\documentclass[12pt,a4paper]{article}
\usepackage{amsmath}
\usepackage{amssymb}
\usepackage{t1enc}
\usepackage[latin1]{inputenc}
\usepackage[english]{babel}
\usepackage{amsfonts}
\usepackage{latexsym}
\usepackage{bm}
\usepackage{setspace}
\usepackage[toc,page]{appendix}
\usepackage{graphicx}
\usepackage{enumerate}
\usepackage[numbers,sort&compress]{natbib}
\usepackage{tikz}

\usepackage{theoremref}
\usepackage{lineno}

\usepackage{mathtools}

\newtheorem{theorem}{Theorem}

\newtheorem{lemma}{Lemma}

\newtheorem{construction}{Construction}
\newtheorem{problem}{Problem}

\allowdisplaybreaks

\begin{document}
\title{Solution to a problem on isolation of cliques in uniform hypergraphs 
\medskip\medskip
}

\author{Peter Borg \\[2mm]
\normalsize Department of Mathematics \\
\normalsize Faculty of Science \\
\normalsize University of Malta\\
\normalsize Malta\\
\normalsize \texttt{peter.borg@um.edu.mt}
}

\date{}
\maketitle

\begin{abstract}
A copy of a hypergraph $F$ is called an $F$-copy. Let $K_k^r$ denote the complete $r$-uniform hypergraph whose vertex set is $[k] = \{1, \dots, k\}$ (that is, the edges of $K_k^r$ are the $r$-element subsets of $[k]$). Given an $r$-uniform $n$-vertex hypergraph $H$, the $K_k^r$-isolation number of $H$, denoted by $\iota(H, K_k^r)$, is the size of a smallest subset $D$ of the vertex set of $H$ such that the closed neighbourhood $N[D]$ of $D$ intersects the vertex sets of the $K_k^r$-copies contained by $H$ (equivalently, $H-N[D]$ contains no $K_k^r$-copy). 
In this note, we show that if $2 \leq r \leq k$ and $H$ is connected, then $\iota(H, K_k^r) \leq \frac{n}{k+1}$ unless $H$ is a $K_k^r$-copy or $k = r = 2$ and $H$ is a $5$-cycle. This solves a recent problem of Li, Zhang and Ye. The result for $r = 2$ (that is, $H$ is a graph) was proved by Fenech, Kaemawichanurat and the author, and is used to prove the result for any $r$. The extremal structures for $r = 2$ were determined by various authors. We use this to determine the extremal structures for any $r$.  
\end{abstract}

\section{Introduction} \label{Introsection}
Unless stated otherwise, we use capital letters such as $X$ to denote sets or graphs, and small letters such as $x$ to denote non-negative integers or elements of a set. The set of positive integers is denoted by $\mathbb{N}$. For $n \geq 1$, $[n]$ denotes the set $\{1, \dots, n\}$ (that is, $\{i \in \mathbb{N} \colon i \leq n\}$). We take $[0]$ to be the empty set $\emptyset$. Arbitrary sets are taken to be finite. A set of sets is called a \emph{family}. A set of size $k$ is called a \emph{$k$-element set} or simply a \emph{$k$-set}. For a set $X$, the \emph{power set of $X$} (the family of subsets of $X$) is denoted by $2^X$, and the family of $k$-element subsets of $X$ is denoted by ${X \choose k}$ (that is, ${X \choose k} = \{A \subseteq X \colon |X| = k\}$). For standard terminology in graph theory, we refer the reader to \cite{West}. Most of the graph terminology used here is defined in \cite{Borg}. 

A \emph{hypergraph} $H$ is a pair $(X, Y)$ such that $X$ is a set denoted by $V(H)$ and called the \emph{vertex set of $H$}, and $Y$ is a subfamily of $2^X$ denoted by $E(H)$ and called the \emph{edge set of $H$}. A member of $V(H)$ is called a \emph{vertex of $H$}, and a member of $E(H)$ is called a \emph{hyperedge of $H$} or simply an \emph{edge of $H$}. If $|V(H)| = n$, then $H$ is said to be an \emph{$n$-vertex hypergraph}. If $E(H) \subseteq {V(H) \choose r}$, then $H$ is said to be \emph{$r$-uniform}. A graph is a $2$-uniform hypergraph. An $r$-uniform hypergraph is also called an \emph{$r$-graph}. We may represent an edge $\{v, w\}$ by $vw$. If $v, w \in e \in E(H)$ and $v \neq w$, then $w$ is called a \emph{neighbour of $v$ in $H$}. If $v \in e \in E(H)$, then $e$ is said to be \emph{incident to $v$ in $H$}. For $v \in V(H)$, the set of neighbours of $v$ in $H$ is denoted by $N_{H}(v)$, and the set $N_{H}(v) \cup \{v\}$ is denoted by $N_{H}[v]$ and called the \emph{closed neighbourhood of $v$ in $H$}. For $X \subseteq V(H)$, the set $\bigcup_{v \in X} N_H[v]$ is denoted by $N_H[X]$ and called the \emph{closed neighbourhood of $X$ in $H$}, the hypergraph $(X, E(H) \cap 2^X)$ is denoted by $H[X]$ and called the \emph{subhypergraph of $H$ induced by $X$}, and the hypergraph $H[V(H) \setminus X]$ is denoted by $H - X$. Where no confusion arises, the subscript $H$ may be omitted; for example, $N_H(v)$ may be abbreviated to $N(v)$. 

If $F$ and $H$ are hypergraphs, $f : V(F) \rightarrow V(H)$ is a bijection, and $E(H) = \{\{f(v) \colon v \in e\} \colon e \in E(F)\}$, then we say that $H$ is a \emph{copy of $F$} or that $H$ is \emph{isomorphic to $F$}, and we write $H \simeq F$. Thus, a copy of $F$ is a hypergraph obtained by relabelling the vertices of $F$. We also call it an \emph{$F$-copy}. If $F$ and $H$ are hypergraphs such that $V(F) \subseteq V(H)$ and $E(F) \subseteq E(H)$, then $F$ is called a \emph{subhypergraph of $H$}, and we say that \emph{$H$ contains $F$}.  

The $r$-graph $([k], {[k] \choose r})$ is denoted by $K_k^r$ and called a \emph{$k$-clique}. For $r = 2$, $K_k^r$ is abbreviated to $K_k$. We call a $K_k^r$-copy contained by an $r$-graph $H$ a \emph{$k$-clique of $H$}. For $n \geq 3$, the graph $([n], \{\{1,2\}, \{2, 3\}, \dots, \{n-1,n\}, \{n,1\}\})$ is denoted by $C_n$. A copy of $C_n$ is called an \emph{$n$-cycle} or simply a \emph{cycle}. A hypergraph $H$ is said to be \emph{connected} if for every $v, w \in V(H)$ with $v \neq w$, there exist some $e_1, \dots, e_t \in E(H)$ such that $v \in e_1$, $w \in e_t$ and $e_i \cap e_{i+1} \neq \emptyset$ for each $i \in [t-1]$. 

If $D \subseteq V(H) = N[D]$, then $D$ is called a \emph{dominating set of $H$}. The size of a smallest dominating set of $H$ is called the \emph{domination number of $H$} and denoted by $\gamma(H)$. If $\mathcal{F}$ is a set of hypergraphs and $F$ is a copy of a hypergraph in $\mathcal{F}$, then we call $F$ an \emph{$\mathcal{F}$-graph}. If $D \subseteq V(H)$ such that $N[D]$ intersects the vertex sets of the $\mathcal{F}$-graphs contained by $H$, then $D$ is called an \emph{$\mathcal{F}$-isolating set of $H$}. Note that $D$ is an $\mathcal{F}$-isolating set of $H$ if and only if $H-N[D]$ contains no $\mathcal{F}$-graph. It is to be assumed that $(\emptyset, \emptyset) \notin \mathcal{F}$. Let $\iota(H, \mathcal{F})$ denote the size of a smallest $\mathcal{F}$-isolating set of $H$. If $\mathcal{F} = \{F\}$, then we may replace $\mathcal{F}$ in these defined terms and notation by $F$. Clearly, for $r \geq 2$, $D$ is a $K_1^r$-isolating set of $H$ if and only if $D$ is a dominating set of $H$, so $\gamma(H) = \iota(H, K_1^r)$. Trivially, $\iota(H, \mathcal{F}) \leq \gamma(H)$.

The study of isolating sets of graphs was introduced by Caro and Hansberg~\cite{CaHa17}. It is a natural generalization of the study of dominating sets \cite{C, CH, HHS, HHS2, HL, HL2}. One of the earliest results in this field is the upper bound $n/2$ of Ore \cite{Ore} on the domination number of any connected $n$-vertex graph $G \not\simeq K_1$ (see \cite{HHS}). While deleting the closed neighbourhood of a dominating set yields the graph with no vertices, deleting the closed neighbourhood of a $K_2$-isolating set yields a graph with no edges. In the literature, a $K_2$-isolating set is also called a \emph{vertex-edge dominating set}. Consider any connected $n$-vertex graph $G$. Caro and Hansberg~\cite{CaHa17} proved that $\iota(G, K_2) \leq n/3$ unless $G \simeq K_2$ or $G \simeq C_5$. This was independently proved by \.{Z}yli\'{n}ski \cite{Z} and solved a problem in \cite{BCHH}. Fenech, Kaemawichanurat and the present author~\cite{BFK} proved the following generalization, which solved a problem in \cite{CaHa17}. 
\begin{theorem}[\cite{BFK}] \label{graphver} If $k \geq 1$ and $G$ is a connected $n$-vertex graph, then, unless either $G \simeq K_k$ or $k = 2$ and $G \simeq C_5$,
\begin{equation} \iota(G, K_k) \leq \frac{n}{k+1}. \label{BFKbound}
\end{equation} 
Moreover, there exists a graph $B_{n,k}$ such that $\iota(B_{n,k}, K_k) = \lfloor n/(k+1) \rfloor$.
\end{theorem}
An explicit construction of $B_{n,k}$ is given in \cite{BFK} and generalized in Construction~\ref{extremalhyp} below. Ore's result is the case $k = 1$, and the result of Caro and Hansberg and of \.{Z}yli\'{n}ski is the case $k = 2$. The graphs attaining the bound in (\ref{BFKbound}) are determined in \cite{FJKR, PX} for $k = 1$, in \cite{BG, LMS} for $k = 2$, in \cite{CCZ} for $k = 3$, and in \cite{CCZ2} for $k \geq 4$. Other isolation bounds of this kind in terms of $n$ are given in \cite{BBS, Borg, Borgcon, Borgrsc, CX, Yan, ZW2, ZW3}. It is worth mentioning that domination and isolation have been particularly investigated for maximal outerplanar graphs \cite{BK, BK2, CaWa13, CaHa17, Ch75, DoHaJo16, DoHaJo17, HeKa18, LeZuZy17, Li16, MaTa96, KaJi, To13}, mostly due to connections with Chv\'{a}tal's Art Gallery Theorem \cite{Ch75}. As in the development of domination, isolation is expanding in various directions, such as total isolation \cite{BGH, CAW} and isolation games~\cite{BDJKR}. 

Li, Zhang and Ye \cite{LZY} asked for a hypergraph version of Theorem~\ref{graphver}. More precisely, they asked for the best possible upper bound on $\iota(H, K_k^r)$ for connected $r$-graphs $H$ \cite[Problems~3.1 and 3.2]{LZY}, and they proved that $\iota(H, K_r^r) \leq n/r$, and asked if $\iota(H, K_r^r) \leq n/(2r-1)$ (unless $H$ is a member of a set of exceptional $r$-graphs). We provide an answer in Theorem~\ref{hypver}. In order to state our results, we need the following construction.

\begin{construction} \label{extremalhyp} \emph{Consider any $n, k, r \in \mathbb{N}$ with $2 \leq r \leq n$, and any connected $k$-vertex $r$-graph $F$. By the division algorithm, there exist $q, s \in \{0\} \cup \mathbb{N}$ such that $n = q(k+1) + s$ and $0 \leq s \leq k$. Let $Q_{n,k}$ be a set of size $q+s$, and let $v_1, \dots, v_{q+s}$ be the elements of $Q_{n,k}$. If $q \geq 1$, then let $F_1, \dots, F_q$ be copies of $F$ such that the $q+1$ sets $V(F_1), \dots, V(F_q)$ and $Q_{n,k}$ are pairwise disjoint, and for each $i \in [q]$, let $\emptyset \neq \mathcal{W}_i \subseteq \{e \in {\{v_i\} \cup V(F_i) \choose r} \colon v_i \in e\}$, and let $H_i$ be the $r$-graph with $V(H_i) = \{v_i\} \cup V(F_i)$ and $E(H_i) = E(F_i) \cup \mathcal{W}_i$. If either $q = 0$ and $H$ is an $n$-vertex $r$-graph that is not an $F$-copy, or $q \geq 1$, $T$ is a connected $r$-graph with $V(T) = Q_{n,k}$, $T'$ is a connected $r$-graph such that $\{v_i \colon i \in [q+s] \setminus [q]\} \subseteq V(T') \subseteq \{v_i \colon i \in [q+s] \setminus [q]\} \cup V(H_q)$ and $v_q \in e$ for each $e \in E(T')$, and $H$ is the $r$-graph with $V(H) = V(T') \cup \bigcup_{i=1}^q V(H_i)$ and $E(H) = E(T) \cup E(T') \cup \bigcup_{i=1}^q E(H_i)$, then we say that $H$ is an \emph{$(n,F)$-good $r$-graph} with \emph{quotient $r$-graph $T$} and \emph{remainder $r$-graph $T'$}, and for each $i \in [q]$, we call $H_i$ an \emph{$F$-constituent of $H$}, and we call $v_i$ the \emph{$F$-connection of $H_i$ in $H$}. We say that an $(n,F)$-good $r$-graph is \emph{pure} if its remainder $r$-graph has no vertices (so $s = 0$). Clearly, an $(n,F)$-good $r$-graph is a connected $n$-vertex $r$-graph.}
\end{construction}

In the next section, we prove the following result.

\begin{theorem} \label{hypver} If $2 \leq r \leq k$ and $H$ is a connected $n$-vertex $r$-graph, then, unless either $H \simeq K_k^r$ or $k = r = 2$ and $H \simeq C_5$,
\begin{equation} \iota(H, K_k^r) \leq \frac{n}{k+1}. \label{hypbound}
\end{equation} 
Moreover, $\iota(H, K_k^r) = \lfloor n/(k+1) \rfloor$ if $H$ is $(n, K_k^r)$-good. 
\end{theorem}
As pointed out above, the graphs attaining the bound in (\ref{BFKbound}) have been completely determined. They are the $r$-graphs attaining the bound in (\ref{hypbound}) for $r = 2$. We determine the $r$-graphs attaining the bound in (\ref{hypbound}) for $r \geq 3$. 

In \cite{CCZ}, Chen, Cui and Zhang defined $10$ connected $8$-vertex graphs $A_1, \dots, A_{10}$ having the same vertex set $\{a_1, \dots, a_8\}$, and proved that the cycle isolation bound $n/4$ in \cite{Borg} is attained by a graph $G \not\simeq K_3$ if and only if $G$ is a pure $(n,K_3)$-good graph or a $\{C_4, A_1, \dots, A_{10}\}$-graph. Consequently, they also proved that the bound in (\ref{BFKbound}) is attained for $k = 3$ if and only if $G$ is a pure $(n,K_3)$-good graph or a $\mathcal{G}_3$-graph, where $\mathcal{G}_3 = \{A_i \colon i \in [10] \setminus \{2\}\}$. In \cite{CCZ2}, Chen, Cui and Zhong treated the case $k \geq 4$. They defined a connected $10$-vertex graph $A$ with vertex set $\{a_1, \dots, a_{10}\}$, and $k+2$ connected $(2k+2)$-vertex graphs $A_k^1, \dots, A_k^{k+2}$ having the same vertex set $\{a_1, \dots, a_{2k+2}\}$. Let $\mathcal{G}_4 = \{A, A_4^1, \dots, A_4^{6}\}$ and $\mathcal{G}_k = \{A_k^1, \dots, A_k^{k+2}\}$ for $k \geq 5$. They proved that for $k \geq 4$, the bound in (\ref{BFKbound}) is attained if and only if $G$ is a pure $(n,K_k)$-good graph or a $\mathcal{G}_k$-graph. Therefore, the results in \cite{CCZ, CCZ2} sum up as follows.

\begin{theorem}[\cite{CCZ, CCZ2}] \label{graphverext} For $k \geq 3$, equality in (\ref{BFKbound}) holds if and only if $G$ is a pure $(n,K_k)$-good graph or a $\mathcal{G}_k$-graph.
\end{theorem}

Let $e_3^1 = \{a_1, a_2, a_3\}$, $e_3^2 = \{a_1, a_2, a_5\}$, $e_3^{3} = \{a_1, a_3, a_5\}$, $e_3^4 = \{a_1, a_5, a_6\}$, $e_3^5 = \{a_1, a_5, a_7\}$, $e_3^6 = \{a_2, a_3, a_4\}$, $e_3^{7} = \{a_2, a_4, a_8\}$, $e_3^8 = \{a_3, a_4, a_8\}$, $e_3^9 = \{a_4, a_7, a_8\}$, $e_3^{10} = \{a_5, a_6, a_7\}$, $e_3^{11} = \{a_5, a_6, a_8\}$, $e_3^{12} = \{a_5, a_7, a_8\}$ and $e_3^{13} = \{a_6, a_7, a_8\}$. Let $\mathcal{E}_3^1 = \{e_3^i \colon i \in \{2, 5, 6, 8, 13\}\}$, $\mathcal{E}_3^2 = \{e_3^i \colon i \in \{1, 4, 6, 9, 10, 13\}\}$, $\mathcal{E}_3^3 = \mathcal{E}_3^2 \cup {\{a_5, a_6, a_7, a_8\} \choose 3}$, $\mathcal{E}_3^4 = \mathcal{E}_3^3 \setminus \{e_3^{10}\}$, $\mathcal{E}_3^5 = \mathcal{E}_3^3 \setminus \{e_3^{11}\}$, $\mathcal{E}_3^6 = \mathcal{E}_3^3 \setminus \{e_3^{12}\}$, $\mathcal{E}_3^7 = \mathcal{E}_3^3 \setminus \{e_3^{13}\}$, $\mathcal{E}_3^8 = \mathcal{E}_3^3 \setminus \{e_3^{10}, e_3^{11}\}$, $\mathcal{E}_3^9 = \mathcal{E}_3^3 \setminus \{e_3^{12}, e_3^{13}\}$ and $\mathcal{E}_3^{10} = \mathcal{E}_3^2 \cup \{e_3^3, e_3^7\}$. For each $i \in [10]$, let $H_3^i = (\{a_1, \dots, a_8\}, \mathcal{E}_3^i)$. Let $\mathcal{H}_3^3 = \{H_3^i \colon i \in [10]\}$. Let $e_4^1 = \{a_1, a_2, a_5, a_6\}$, $e_4^2 = \{a_1, a_6, a_7, a_{10}\}$, $e_4^3 = \{a_2, a_3, a_4, a_5\}$, $e_4^4 = \{a_3, a_4, a_8, a_9\}$ and $e_4^5 = \{a_7, a_8, a_9, a_{10}\}$. Let $H_4^4 = (\{a_1, \dots, a_{10}\}, \{e_4^1, \dots, e_4^5\})$ and $H_4^3 = (\{a_1, \dots, a_{10}\}, {e_4^1 \choose 3} \cup \dots \cup {e_4^5 \choose 3})$. Let $\mathcal{H}_4^3 = \{H_4^3\}$ and $\mathcal{H}_4^4 = \{H_4^4\}$. In the next section, we also prove the following result.

\begin{theorem} \label{hypverext} For $3 \leq r \leq k$, equality in (\ref{hypbound}) holds if and only if $H$ is a pure $(n,K_k^r)$-good $r$-graph or $3 \leq k \leq 4$ and $H$ is an $\mathcal{H}_k^r$-graph.
\end{theorem}
We convert the $r$-graph setting to a graph setting. This enables us to obtain Theorem~\ref{hypver} from Theorem~\ref{graphver}, and to obtain Theorem~\ref{hypverext} from Theorem~\ref{graphverext}. 

\section{Proofs} \label{Proofsection}

We now start working towards proving Theorems~\ref{hypver} and \ref{hypverext}.

For a family $\mathcal{A}$ of sets, the family $\bigcup_{A \in \mathcal{A}} {A \choose s}$ is denoted by $\partial_s (\mathcal{A})$ and called the \emph{$s$th shadow of $\mathcal{A}$}. For a hypergraph $H$, we denote by $H^{(s)}$ the $s$-graph with vertex set $V(H)$ and edge set $\partial_s (E(H))$.

\begin{lemma} \label{keylemma} Let $2 \leq s \leq r \leq k$ and let $H$ be an $r$-graph. \\
(i) For any $D \subseteq V(H)$, $N_H[D] = N_{H^{(s)}}[D]$. \\
(ii) If $D$ is a $K_k^s$-isolating set of $H^{(s)}$, then $D$ is a $K_k^r$-isolating set of $H$.\\
(iii) $E(H) \subseteq \{V(R) \colon R \mbox{ is an $r$-clique of } H^{(s)}\}$. \\
(iv) If $e \in E(H^{(s)})$ and $H^{(s)}$ contains only one $K_r^s$-copy $F$ with $e \in E(F)$, then $V(F) \in E(H)$.\\
(v) If $e \in E(H)$ and $H$ has no $k$-clique $F$ with $e \in E(F)$, then there exists no $k$-graph $I$ with $I^{(r)} = H$.
\end{lemma}
\textbf{Proof.} Let $D \subseteq V(H)$. We have $D \subseteq N_H[D] \cap N_{H^{(s)}}[D]$. Let $v \in V(H)$. Suppose $v \in N_H[D] \setminus D$. Then, $v \in N_H[u]$ for some $u \in D$, so $u, v \in e$ for some $e \in E(H)$. Let $e' \subseteq e$ such that $u, v \in e'$ and $|e| = s$. Then, $e' \in H^{(s)}$, so $v \in N_{H^{(s)}}[u]$. Thus, $N_H[D] \subseteq N_{H^{(s)}}[D]$. Now suppose $v \in N_{H^{(s)}}[D] \setminus D$. Then, $u, v \in e$ for some $u \in D$ and $e \in E(H^{(s)})$. Since $e \subseteq e'$ for some $e' \in E(H)$, $v \in N_H[u]$. Thus, $N_{H^{(s)}}[D] \subseteq N_{H}[D]$. Since $N_H[D] \subseteq N_{H^{(s)}}[D]$, (i) follows.

Suppose that $D$ is a $K_k^s$-isolating set of $H^{(s)}$ and that $H$ contains a copy $B$ of $K_k^r$. Then, $B^{(s)}$ is a copy of $K_k^s$ contained by $H^{(s)}$. Thus, $N_{H^{(s)}}[D] \cap V(B^{(s)}) \neq \emptyset$. By (i), $N_{H}[D] \cap V(B^{(s)}) \neq \emptyset$. Since $V(B^{(s)}) = V(B)$, (ii) follows.

If $e \in E(H)$, then $(e, {e \choose s})$ is an $r$-clique of $H^{(s)}$. This yields (iii).

Suppose that $e \in E(H^{(s)})$ and $H^{(s)}$ contains only one $K_r^s$-copy $F$ with $e \in E(F)$. We have $e \subseteq e'$ for some $e' \in E(H)$. Let $F' = (e', {e' \choose s})$. Then, $e \in E(F')$ and $F'$ is a $K_r^s$-copy contained by $H^{(s)}$, so $F' = F$. We have $V(F) = V(F') = e' \in E(H)$, so (iv) is proved.

Suppose that $e \in E(H)$ and $I$ is a $k$-graph with $I^{(r)} = H$. Then, $e \subseteq e'$ for some $e' \in E(I)$. Let $F' = (e', {e' \choose r})$. Then, $e \in E(F')$ and $F'$ is a $K_k^r$-copy contained by $H$. This yields (v).~\hfill{$\Box$}
\\

The converse of Lemma~\ref{keylemma} (ii) is false. Indeed, if $s < r < k$ and $H = ([k], {[k] \choose r} \setminus \{[r]\})$, then $H$ contains no $K_k^r$-copy and $H^{(s)}$ is a $K_k^s$-copy (so $\emptyset$ is a $K_k^r$-isolating set of $H$ but not a $K_k^s$-isolating set of $H^{(s)}$).
\\

\noindent
\textbf{Proof of Theorem~\ref{hypver}.} If $r = 2$, then the result is given by Theorem~\ref{graphver}. Suppose $r \geq 3$. If $n \leq k$, then $\iota(H, K_k^r) = 0$ unless $H \simeq K_k^r$. Suppose $n \geq k + 1$. Let $G$ be the graph $H^{(2)}$. Since $H$ is connected, $G$ is connected. Since $n \geq k + 1$, $G \not\simeq K_k$. Since $r \geq 3$, $G \not\simeq C_5$. Let $D$ be a smallest $K_k$-isolating set of $G$. By Theorem~\ref{graphver}, $|D| \leq n/(k+1)$. By Lemma~\ref{keylemma} (ii), $D$ is a $K_k^r$-isolating set of $H$. This yields (\ref{hypbound}). 

Now suppose that $H$ is an $(n,K_k^r)$-good $r$-graph with exactly $q$ $K_k^r$-constituents as in Construction~\ref{extremalhyp}. Then, $q = \lfloor n/(k+1) \rfloor$. If $q = 0$, then $\iota(H,K_k^r) = 0$. Suppose $q \geq 1$. Then, $\{v_1, \dots, v_q\}$ is a $K_k^r$-isolating set of $H$. If $D$ is a $K_k^r$-isolating set of $H$, then, since $H_1 - v_1, \dots, H_q - v_q$ are copies of $K_k^r$, we have $D \cap V(H_i) \neq \emptyset$ for each $i \in [q-1]$, and $D \cap (V(H_q) \cup V(T')) \neq \emptyset$. Therefore, $\iota(H, K_k^r) = q$.~\hfill{$\Box$}
\\

\noindent
\textbf{Proof of Theorem~\ref{hypverext}.} We first settle the necessary condition. Thus, suppose that $H$ attains the bound in (\ref{hypbound}). Let $G$ and $D$ be as in the proof of Theorem~\ref{hypver}. By Theorem~\ref{graphver}, $|D| \leq n/(k+1)$. By Lemma~\ref{keylemma}~(ii), $D$ is a $K_k^r$-isolating set of $H$, so $|D| \geq \iota(H, K_k^r)$. We have $n/(k+1) = \iota(H, K_k^r) \leq |D| \leq n/(k+1)$, so $|D| = n/(k+1)$. By Theorem~\ref{graphverext}, $G$ is a pure $(n,K_k)$-good graph or a $\mathcal{G}_k$-graph. 

Suppose that $G$ is a pure $(n,K_k)$-good graph. We may assume that $G$ is as in Construction~\ref{extremalhyp} (with $F = K_k$). Thus, $Q_{n,k} = \{v_1, \dots, v_q\}$, $V(G) = Q_{n,k} \cup \bigcup_{i=1}^q V(F_i)$, and for each $i \in [q]$, we have $V(F_i) \simeq K_k$, $N_G[V(F_i)] = V(F_i) \cup \{v_i\}$, and hence $N_H[V(F_i)] = V(F_i) \cup \{v_i\}$ by Lemma~\ref{keylemma}~(i). Let $Q = Q_{n,k}$. Suppose $H[V(F_j)] \not\simeq K_k^r$ for some $j \in [q]$. Suppose $q \geq 2$. Let $Q' = Q \setminus \{j\}$. Since $G[Q]$ is connected, $H[Q]$ is connected, so $v_j \in N_H[Q']$. We obtain that $Q'$ is a $K_k^r$-isolating of $H$. We have $|Q'| < q = n/(k-1)$, contradicting $\iota(H, K_k^r) = n/(k+1)$. Thus, $H[V(F_i)] \simeq K_k^r$ for each $i \in [q]$, and hence $H$ is a pure $(n,K_k^r)$-good $r$-graph. Now suppose $q = 1$. We have $1 = q = n/(k+1) = \iota(H,K_k^r)$, so $H$ contains a $K_k^r$-copy $I$. Since $n = k + 1 = |V(I)| + 1$, $V(H) = V(I) \cup \{v\}$ for some $v \in V(H) \setminus V(I)$. Since $H$ is connected, $H$ is a pure $(n,K_k^r)$-good $r$-graph.

Now suppose that $G$ is a $\mathcal{G}_k$-graph. We may assume that $G \in \mathcal{G}_k$. Suppose $k \geq 4$. Let $J \in \{A_k^1, \dots, A_k^{k+2}\}$. Then, $a_{k+2} \in N_J[a_1] \subseteq \{a_1, \dots, a_{k+2}\}$ and $a_1 \in N_J[a_{k+2}] \subseteq \{a_1, a_{k+2}, \dots, a_{2k+2}\}$ (see \cite{CCZ2}). Thus, $a_1a_{k+2} \in E(J)$ and $J[\{v, a_1, a_{k+2}\}] \not\simeq K_3$ for each $v \in V(J) \setminus \{a_1, a_{k+2}\}$. Since $r \geq 3$, $J$ contains no $K_r$-copy $F$ with $a_1a_{k+2} \in E(F)$. By Lemma~\ref{keylemma} (v), $H^{(2)} \neq J$, so $G \neq J$. Therefore, $k = 4$ and $G = A$. The $4$-cliques of $G$ are $G[e_4^1], \dots, G[e_4^5]$, and the set of $3$-cliques of $G$ is $\bigcup_{i=1}^5 \{G[T] \colon T \in {e_4^i \choose 3}\}$ (see \cite{CCZ2}). By Lemma~\ref{keylemma}~(iii), $E(H) \subseteq E(H_4^r)$. Let $a_1' = a_8$, $a_2' = a_3$, $a_3' = a_7$, $a_4' = a_{1}$ and $a_5' = a_2$. By Lemma~\ref{keylemma} (i), $N_H[v] = N_G[v]$ for each $v \in V(H)$. For each $i \in [5]$, $H - N_H[a_i'] = H - N_G[a_i'] = H[e_4^i]$, so $H[e_4^i] \simeq K_4^r$ as $\iota(H,K_4^r) = n/(k+1) = 10/5 = 2$. Therefore, $\bigcup_{i = 1}^5 {e_4^i \choose r} \subseteq E(H)$, and hence $H = H_4^r$.

Now suppose $k = 3$. Since $3 \leq r \leq k$, $r = 3$. Let $J \in \{A_6, A_7, A_8, A_9, A_{10}\}$. Then, $a_{5} \in N_J[a_1] \subseteq \{a_1, \dots, a_5\}$ and $a_1 \in N_J[a_5] \subseteq \{a_1, a_5, \dots, a_{8}\}$ (see \cite{CCZ}). Thus, $a_1a_5 \in E(J)$ and $J[\{v, a_1, a_5\}] \not\simeq K_3$ for each $v \in V(J) \setminus \{a_1, a_{5}\}$. By Lemma~\ref{keylemma} (v), $H^{(2)} \neq J$, so $G \neq J$. Thus, $G = A_j$ for some $j \in \{1, 3, 4, 5\}$. Let $X = \{1, 3, 4, 5\}$. For each $i \in X$, let $\mathcal{K}_i$ be the family of vertex sets of the $3$-cliques of $A_i$, and let $S_i = \{(e, V(F)) \colon F \mbox{ is the only $3$-clique of } A_i \mbox{ with } e \in E(F)\}$. Let 
\begin{align} S_5' = \{&(a_1a_2, e_3^1), (a_3a_5, e_3^3), (a_1a_6, e_3^4), (a_3a_4, e_3^6), (a_2a_8, e_3^7), \nonumber \\
&(a_4a_7, e_3^9), (a_5a_7, e_3^{10}), (a_6a_8, e_3^{13})\}. \nonumber
\end{align}
It can be checked that $S_5' \subseteq S_5$ and $\mathcal{K}_5 = \{T \colon (e, T) \in S_5' \mbox{ for some } e \in E(A_5)\}$. Thus, if $j = 5$, then by Lemma~\ref{keylemma} (iii) and (iv), $E(H) = \mathcal{K}_5$, and hence $H = H_3^{10}$. Since $E(A_3) \subseteq E(A_5)$, we similarly obtain $H = H_3^2$ if $j = 3$. Let $S_1' = \{(a_2a_5, e_3^2), (a_1a_7, e_3^5), (a_2a_3, e_3^6), (a_3a_8, e_3^8), (a_6a_7, e_3^{13})\}$. Since $S_1' \subseteq S_1$ and $\mathcal{K}_1 = \{T \colon$ $(e, T) \in S_1' \mbox{ for some } e \in E(A_1)\}$, we obtain $H = H_3^1$ if $j = 1$. Finally, suppose $j = 4$. Let $S_4' = \{(a_1a_2, e_3^1), (a_1a_6, e_3^4)$, $(a_3a_4, e_3^6), (a_4a_7, e_3^9)\}$. Since $S_4' \subseteq S_4$, $e_3^1, e_3^4, e_3^6, e_3^9$ are hyperedges of $H$ by Lemma~\ref{keylemma}~(iv). Let $\mathcal{E}^*$ be the set of these $4$ hyperedges, and let $Z = \{a_5, a_6, a_7, a_8\}$ and $\mathcal{E}' = {Z \choose 3}$. We have $\mathcal{E}^* \subseteq E(H)$, $\mathcal{E}_3^3 = \mathcal{E}^* \cup \mathcal{E}'$ and $\mathcal{K}_4 = \mathcal{E}_3^3$. By Lemma~\ref{keylemma}~(iii), $E(H) \subseteq \mathcal{E}_3^3$. Since $a_5a_8 \in E(G)$, we have $e_3^{11} \in E(H)$ or $e_3^{12} \in E(H)$. Suppose $e_3^{11} \in E(H)$. Since $a_6a_7 \in E(G)$, we have $e_3^{10} \in E(H)$ or $e_3^{13} \in E(H)$. If $e_3^{10} \in E(H)$, then $H \in \{H_3^3, H_3^6, H_3^7, H_3^9\}$. If $e_3^{13} \in E(H)$, then since $a_5a_7 \in E(G)$, we have $e_3^{10} \in E(H)$ or $e_3^{12} \in E(H)$, so $H \in \{H_3^3, H_3^4, H_3^6\}$. Now suppose $e_3^{11} \notin E(H)$. Then, $e_3^{12} \in E(H)$. Since $e_3^{11} \notin E(H)$ and $a_6a_8 \in E(G)$, $e_3^{13} \in E(H)$. Thus, $H \in \{H_3^5, H_3^8\}$. 

We now settle the sufficient condition. By Theorem~\ref{hypver}, $\iota(H, K_k^r) = n/(k+1)$ if $H$ is a pure $(n,K_k^r)$-good $r$-graph. Now suppose $3 \leq k \leq 4$ and $H \in \mathcal{H}_k^r$. 
It is easily checked that if $3 = r = k$, then $H - N_H[a_i]$ contains a $K_3^3$-copy for each $i \in [8]$, so we have $1 < \iota(H, K_3^3) \leq n/(k+1) = 2$, and hence $\iota(H, K_3^3) = n/(k+1)$. Similarly, if $3 \leq r \leq k = 4$, then $H - N_H[a_i]$ contains a $K_k^r$-copy for each $i \in [10]$, so $\iota(H, K_k^r) = 2 = n/(k+1)$.~\hfill{$\Box$}

\section{The case $k < r$}

The problem of obtaining best possible upper bounds on $\iota(H, K_k^r)$ is fundamentally different for $k < r$. In this case, $K_k^r$ has no edges, and hence if $k \geq 2$, then $K_k^r$ is not connected. In general, given a set $\mathcal{F}$ of hypergraphs, certain desirable properties of $\mathcal{F}$-isolating sets are not guaranteed if some members of $\mathcal{F}$ are not connected. In particular, if $\mathcal{H}$ is the set of components of $H$, then $\iota(H, \mathcal{F}) = \sum_{I \in \mathcal{H}} \iota(I, \mathcal{F})$ if the members of $\mathcal{F}$ are connected, but $\iota(H, \mathcal{F})$ may not be $\sum_{I \in \mathcal{H}} \iota(I, \mathcal{F})$ otherwise; see \cite[Section~2]{Borgrsc}.  

We pose the following problem.

\begin{problem}\label{problem} For $r \geq 3$ and $1 \leq k < r \leq n$, what is the smallest rational number $c = c(n,k,r)$ such that $\iota(H, K_k^r) \leq c n$ for every connected $n$-vertex $r$-graph $H$?
\end{problem}

As pointed out in Section~\ref{Introsection}, for $k = 1$, Problem~\ref{problem} is the famous domination problem for $r$-graphs. For $r \in \{3, 4\}$, it is shown in \cite{BHT, HS} that $\gamma(H) \leq n/r$, and that this bound is sharp. For $r = 5$, it is shown in \cite{BHT} that $\gamma(H) \leq 2n/9$. 

Problem~\ref{problem} has the following relation with the domination problem. 

\begin{theorem} If $1 \leq k < r$ and $H$ is an $r$-graph, then
\begin{equation} \gamma(H) - k + 1 \leq \iota(H,K_k^r) \leq \gamma(H). \label{k<r relation}
\end{equation}
Moreover, for every $q \geq 1$, there exist two connected $r$-graphs $H$ and $I$ such that $\iota(H,K_k^r) = \gamma(H) = q = \iota(I,K_k^r) = \gamma(I) - k + 1$. 
\end{theorem}
\textbf{Proof.} As pointed out in Section~\ref{Introsection}, $\iota(H,K_k^r) \leq \gamma(H)$ trivially. Since $k < r$, a subset $D$ of $V(H)$ is a $K_k^r$-isolating set of $H$ if and only if $|V(H) \setminus N[D]| \leq k-1$. Let $D$ be a smallest $K_k^r$-isolating set of $H$, and let $D' = V(H) \setminus N[D]$. Then, $|D'| \leq k-1$, $D \cup D'$ is a dominating set of $H$, and hence $\gamma(H) \leq |D \cup D'| = |D| + |D'| \leq \iota(H,K_k^r) + k-1$. Therefore, (\ref{k<r relation}) is proved. 

Let $q \geq 1$ and $n = q(r+1)$. Suppose that $H$ is a pure $(n,K_r^r)$-good $r$-graph (thus having exactly $q$ $K_r^r$-constituents) as in Construction~\ref{extremalhyp} with $\mathcal{W}_i = \{e \in {\{v_i\} \cup V(F_i) \choose r} \colon v_i \in e\}$ for each $i \in [q]$. Let $X = \{v_1, \dots, v_q\}$. Then, $X$ is a dominating set of $H$. If $D_H$ is a $K_k^r$-isolating set of $H$, then since $H_1 - v_1, \dots, H_q - v_q$ are copies of $K_r^r$ (and hence contain copies of $K_k^r$), we have $D_H \cap V(H_i) \neq \emptyset$ for each $i \in [q]$. Thus, we have $q \leq \iota(H,K_k^r) \leq \gamma(H) \leq |X| = q$, and hence $\iota(H,K_k^r) = \gamma(H) = q$. Let $R_1', \dots, R_{k-1}', S_1', \dots, S_{k-1}'$ be pairwise disjoint sets such that for each $i \in [k-1]$, $|R_i'| = r-1$, $|S_i'| = 1$ and $R_i' \cap V(H) = \emptyset = S_i' \cap V(H)$. For each $i \in [k-1]$, let $R_i = \{v_q\} \cup R_i'$ and $S_i = R_i' \cup S_i'$. Let $I$ be the connected ($n + (k-1)r$)-vertex $r$-graph with vertex set $V(H) \cup \bigcup_{i=1}^{k-1} S_i$ and edge set $E(H) \cup \bigcup_{i=1}^{k-1} \{R_i, S_i\}$. We have $V(I) \setminus N_I[X] = \bigcup_{i=1}^{k-1} S_i'$, so $|V(I) \setminus N_I[X]| \leq k-1$, and hence $X$ is a $K_k^r$-isolating set of $I$. As above, if $D_I$ is a dominating set of $I$ or a $K_k^r$-isolating set of $I$, then $D_I$ intersects each of $V(H_1), \dots, V(H_q)$. Thus, $\iota(I,K_k^r) = |X| = q$. Let $D_I$ be a smallest dominating set of $I$. For each $i \in [k-1]$, we have $S_i' \subseteq N_I[D_I]$, so $D_I \cap S_i \neq \emptyset$. Thus, $|D_I| \geq q+k-1$. Since $X \cup \bigcup_{i=1}^{k-1}S_i'$ is a dominating set of $I$, $\gamma(I) = q+k-1 = \iota(I,K_k^r) + k - 1$.~\hfill{$\Box$}


\footnotesize
\begin{thebibliography}{}

\bibitem{BBS} K. Bartolo, P. Borg, D. Scicluna, Isolation of squares in graphs, Discrete Mathematics 347 (2024), paper 114161.

\bibitem{Borg} P. Borg, Isolation of cycles, Graphs and Combinatorics 36 (2020), 631--637.

\bibitem{Borgcon} P. Borg, Isolation of connected graphs, Discrete Applied Mathematics 339 (2023), 154--165.

\bibitem{Borgrsc} P. Borg, Isolation of regular graphs, stars and $k$-chromatic graphs, Discrete Mathematics 349 (2026), paper 114706.


\bibitem{BFK} P. Borg, K. Fenech, P. Kaemawichanurat, Isolation of $k$-cliques, Discrete Mathematics 343 (2020), paper 111879.


\bibitem{BK} P. Borg, P. Kaemawichanurat, Partial domination of maximal outerplanar graphs, Discrete Applied Mathematics 283 (2020), 306--314.

\bibitem{BK2} P. Borg, P. Kaemawichanurat, Extensions of the Art Gallery Theorem, Annals of Combinatorics 27 (2023) 31--50.

\bibitem{BCHH} R. Boutrig, M. Chellali, T.W. Haynes, S.T. Hedetniemi, Vertex-edge domination in graphs, Aequationes Mathematicae 90 (2016), 355--366.

\bibitem{BG} G. Boyer, W. Goddard, Disjoint isolating sets and graphs with maximum isolation number, Discrete Applied Mathematics 356 (2024), 110--116.

\bibitem{BGH} G. Boyer, W. Goddard, M.A. Henning, On total isolation in graphs, Aequationes Mathematicae 99 (2025), 623--633.

\bibitem{BDJKR} B. Bre\v{s}ar, T. Dravec, D.P. Johnston, K. Kuenzel, D.F. Rall, Isolation game on graphs, arXiv:2409.14180 [math.CO].

\bibitem{BHT} C. Bujt\'{a}s, M.A. Henning, Z. Tuza, Transversals and domination in uniform hypergraphs, European Journal of Combinatorics 33 (2012), 62--71.

\bibitem{CaWa13} C.N. Campos, Y. Wakabayashi, On dominating sets of maximal outerplanar graphs, Discrete Applied Mathematics 161 (2013), 330--335.

\bibitem{CAW} Y. Cao, X. An, B. Wu, Total isolation of $k$-cliques in a graph, Discrete Mathematics 348 (2025), paper 114689.

\bibitem{CaHa17} Y. Caro, A. Hansberg, Partial domination - the isolation number of a graph, Filomat 31 (2017), 3925--3944.

\bibitem{CCZ} S. Chen, Q. Cui, J. Zhang, A characterization of graphs with maximum cycle isolation number, Discrete Applied Mathematics 366 (2025), 161--175.

\bibitem{CCZ2} S. Chen, Q. Cui, L. Zhong, A characterization of graphs with maximum $k$-clique isolation number, Discrete Mathematics 348 (2025), 114531.

\bibitem{CX} J. Chen, S. Xu, $P_5$-isolation in graphs, Discrete Applied Mathematics 340 (2023), 331--349.

\bibitem{Ch75} V. Chv\'{a}tal, A combinatorial theorem in plane geometry, Journal of Combinatorial Theory Series~B 18 (1975), 39--41.

\bibitem{C} E.J. Cockayne, Domination of undirected graphs -- A survey, in: Lecture Notes in Mathematics, Volume 642, Springer, 1978, 141--147.

\bibitem{CH} E.J. Cockayne, S.T. Hedetniemi, Towards a theory of domination in graphs, Networks 7 (1977), 247--261.

\bibitem{DoHaJo16} M. Dorfling, J.H. Hattingh, E. Jonck, Total domination in maximal outerplanar graphs II, Discrete Mathematics 339 (2016), 1180--1188.

\bibitem{DoHaJo17} M. Dorfling, J.H. Hattingh, E. Jonck, Total domination in maximal outerplanar graphs, Discrete Applied Mathematics 217 (2017), 506--511.

\bibitem{FJKR} J.F. Fink, M.S. Jacobson, L.F. Kinch, J. Roberts, On graphs having domination number half their order, Periodica Mathematica Hungarica 16 (1985), 287--293.

\bibitem{HHS} T.W. Haynes, S.T. Hedetniemi, P.J. Slater, Fundamentals of Domination in Graphs, Marcel Dekker, Inc., New York, 1998.

\bibitem{HHS2} T.W. Haynes, S.T. Hedetniemi, P.J. Slater (Editors), Domination in Graphs: Advanced Topics, Marcel Dekker, Inc., New York, 1998.

\bibitem{HL2} S.T. Hedetniemi, R.C. Laskar, Bibliography on domination in graphs and some basic definitions of domination parameters, Discrete Mathematics 86 (1990), 257--277.

\bibitem{HL} S.T. Hedetniemi, R.C. Laskar (Editors), Topics on Domination, in: Annals of Discrete Mathematics, Volume 48, North-Holland Publishing Co., Amsterdam, 1991, Reprint of Discrete Mathematics 86 (1990). 

\bibitem{HeKa18} M.A. Henning, P. Kaemawichanurat, Semipaired domination in maximal outerplanar graphs, Journal of Combinatorial Optimization 38 (2019), 911--926.

\bibitem{HS} M.A. Henning, H.C. Swart, Bounds on a generalized domination parameter, Quaestiones Mathematicae 13 (1990), 237--257.

\bibitem{LMS} M. Lema\'{n}ska,  M. Mora, M.J. Souto-Salorio, Graphs with isolation number equal to one third of the order, Discrete Mathematics 347 (2024), paper 113903.

\bibitem{LeZuZy17} M. Lema\'{n}ska, R. Zuazua, P. \.{Z}yli\'{n}ski, Total dominating sets in maximal outerplanar graphs, Graphs and Combinatorics 33 (2017), 991--998.

\bibitem{LZY} M. Li, S. Zhang, C. Ye, Partial domination of hypergraphs, Graphs and Combinatorics 40 (2024), paper 109 .

\bibitem{Li16} Z. Li, E. Zhu, Z. Shao, J. Xu, On dominating sets of maximal outerplanar and planar graphs, Discrete Applied Mathematics 198 (2016), 164--169.

\bibitem{MaTa96} L.R. Matheson, R.E. Tarjan, Dominating sets in planar graphs, European Journal of Combinatorics 17 (1996), 565--568.

\bibitem{Ore} O. Ore, Theory of graphs, in: American Mathematical Society Colloquium Publications, Volume 38, American Mathematical Society, Providence, R.I., 1962.

\bibitem{PX} C. Payan, N.H. Xuong, Domination-balanced graphs, Journal of Graph Theory 6 (1982), 23--32.

\bibitem{To13} S. Tokunaga, Dominating sets of maximal outerplanar graphs, Discrete Applied Mathematics 161 (2013), 3097--3099.

\bibitem{KaJi} S. Tokunaga, T. Jiarasuksakun, P. Kaemawichanurat, Isolation number of maximal outerplanar graphs, Discrete Applied Mathematics 267 (2019), 215--218.

\bibitem{West} D.B. West, Introduction to Graph Theory, second edition, Prentice Hall, 2001. 

\bibitem{Yan} J. Yan, Isolation of the diamond graph, Bulletin of the Malaysian Mathematical Sciences Society 45 (2022), 1169--1181.

\bibitem{ZW2} G. Zhang, B. Wu, $K_{1,2}$-isolation in graphs, Discrete Applied Mathematics 304 (2021), 365--374.

\bibitem{ZW3} G. Zhang, B. Wu, On the cycle isolation number of triangle-free graphs, Discrete Mathematics 347 (2024), paper 114190.

\bibitem{Z} P. \.{Z}yli\'{n}ski, Vertex-edge domination in graphs, Aequationes Mathematicae 93 (2019), 735--742.

\end{thebibliography}
\end{document}